\documentclass[12pt]{amsart}

\usepackage{amssymb,amscd}
\usepackage{verbatim}
\usepackage{graphicx}

\newtheorem*{Whitney towers}{Theorem~\ref{Whitney towers}}
\newtheorem*{h-towers}{Theorems ~\ref{half} \& \ref{$(n)$-solvable}}

\newtheorem*{surgery curves}{Theorem~\ref{surgery curves}}
\newtheorem*{cg=0}{Theorem~\ref{vanish}}

\newtheorem{thm}{Theorem}[section]

\newtheorem{cla}[thm]{Claim}
\theoremstyle{definition}
\newtheorem{defn}[thm]{Definition}

\newtheorem{note}[thm]{Note}

\newtheorem{prob}[thm]{Problem}

\numberwithin{equation}{section}
\numberwithin{figure}{section}

\newcommand{\x}{\times}
\newcommand{\np}{\newpage}
\newcommand{\Z}{\mathbb{Z}}
\newcommand{\N}{\mathbb{N}}

\newcommand{\Q}{\mathbb{Q}}

\newcommand{\R}{\mathbb{R}}

\newcommand{\MM}{\mathcal{M}}

\newcommand{\nl}{\newline}

\title{
A new obstruction for ribbon-moves of 2-knots:  
2-knots fibred by the punctured 3-tori 
and  
2-knots bounded by homology spheres 
}
\author{Eiji Ogasa}
\address{
Computer Science, Meiji Gakuin University, Yokohama, Kanagawa, 244-8539, 
Japan 
}
\email{pqr100pqr100@yahoo.co.jp, ogasa@mail1.meijigakuin.ac.jp}

\thanks{
{\bf Keyword:}
2-knots, 
ribbon-moves of 2-knots, 
$o$-sets, 
the $o$-invariant, 
the $\infty$-chain homology groups, 
the Alexander module, 
the cyclic covering space of the complement of any 2-knot, 
2-knots whose Seifert hypersurfaces are punctured homology spheres, 
fibered 2-knots whose Seifert hypersurfaces are the punctured $T^3$}

\begin{document}

\begin{abstract} 
This paper gives a new obstruction for ribbon-move equivalence of 2-knots.

Let $K$ and $K'$ be 2-knots. 
Let $K$ and $K'$ are ribbon-move equivalent.
One corollary to our main theorem is as follows. 
A 2-dimensional fibered knot 
whose fiber is the punctured 3-dimensional torus 
 is not ribbon-move equivalent to 
any 2-dimensional knot whose Seifert hypersurface is 
a punctured homology sphere. 

\end{abstract}

\maketitle

\section{Introduction}\label{introduction}



\noindent 
In this paper we give a new obstruction for ribbon move-equivalence of 2-knots. 
The author's papers \cite{O1,O2} gave some obstructions. 
This paper gives new results (\S \ref{main}).

One of the new results is as follows.  
This theorem is deduced from our main theorem (Theorem \ref{daiji}). 

\vskip3mm
\begin{thm}\label{saisho} \rm{ (same as Theorem \ref{cor}.)} 
Let $K$ be a 2-dimensional fibered knot 
whose fiber is the punctured 3-dimensional torus. 
Let $P$ be a 2-dimensional knot whose Seifert hypersurface is 
a punctured homology sphere. 
Then $K$ is not ribbon-move equivalent to $P$. 
\end{thm}

\vskip3mm
\noindent{\bf Note.} 
Let $M$ be a closed $n$-manifold. 
A {\it punctured} ({\it manifold}) $M$ is 
a manifold with boundary,

\noindent 
$M-($ an open $n$-ball embedded trivially in $M$).

\vskip3mm

This paper is based on the author's preprint \cite{Op}. 

We work in the smooth category. 


\section{Review of Ribbon-moves of 2-knots}\label{ribbon}


\noindent 
In this section we review the definition of ribbon-moves.  

An   ({\it oriented})  2-({\it dimensional})  {\it knot}   
is a smooth oriented submanifold $K$ of $S^4$ 
which is diffeomorphic to the $2$-sphere. 
We say that 2-knots $K_1$ and $K_2$ are {\it equivalent} 
if there exists an orientation preserving diffeomorphism  
$f:$ $S^4$ $\rightarrow$ $S^4$ 
such that $f(K_1)$=$K_2$  and 
that $f | _{K_1}:$ $K_1$ $\rightarrow$ $K_2$ is 
an orientation preserving diffeomorphism.    
Let $id:S^4$ $\rightarrow$ $S^4$ be the identity. 
We say that 2-knots $K_1$ and $K_2$ are {\it identical}  
if  $id(K_1)$=$K_2$ and $id | _{K_1}:K_1\rightarrow K_2$ is 
an orientation preserving diffeomorphism.

\vskip3mm
If $W$ is a subset of a manifold $Z$, then let 
$\overline{W}$ mean the closure of $W$ in $Z$. 
In this paper we often omit to explain what 
$Z$ is if it is easy to understand what 
$Z$ is. 
\vskip3mm

\vskip3mm
\begin{defn}\label{ribbonmove}     
Let $K_1$ and $K_2$ be 2-knots in $S^4$. 
We say that $K_2$ is obtained from $K_1$ by one {\it ribbon-move } 
if there is a 4-ball $B$ embedded trivially in $S^4$ with the following properties.  

\noindent (1) 
$K_1$ coincides with $K_2$ in $\overline{S^4-B}$. 
This identity map from $\overline{K_1-B}$ to $\overline{K_2-B}$  
is orientation preserving.

\noindent (2) 
$B\cap K_1$ is drawn as in Figure \ref{ribbonmove}.1.    
$B\cap K_2$ is drawn as in Figure \ref{ribbonmove}.2.    

\np
\vskip3mm
\unitlength 0.1in
\begin{picture}(56.10,51.00)(8.50,-51.60)
%
\special{pn 8}%
\special{ar 3510 320 560 250  0.0000000 6.2831853}%
%
\special{pn 20}%
\special{pa 3660 3950}%
\special{pa 3660 2730}%
\special{fp}%
%
\special{pn 20}%
\special{pa 3340 3970}%
\special{pa 3340 2740}%
\special{fp}%
%
\special{pn 20}%
\special{pa 3350 4000}%
\special{pa 3360 3970}%
\special{pa 3384 3949}%
\special{pa 3413 3935}%
\special{pa 3444 3926}%
\special{pa 3475 3921}%
\special{pa 3507 3920}%
\special{pa 3539 3923}%
\special{pa 3570 3930}%
\special{pa 3601 3940}%
\special{pa 3627 3958}%
\special{pa 3647 3983}%
\special{pa 3648 4014}%
\special{pa 3630 4040}%
\special{pa 3603 4058}%
\special{pa 3573 4070}%
\special{pa 3542 4077}%
\special{pa 3510 4080}%
\special{pa 3478 4079}%
\special{pa 3447 4075}%
\special{pa 3416 4066}%
\special{pa 3387 4053}%
\special{pa 3362 4032}%
\special{pa 3350 4003}%
\special{pa 3350 4000}%
\special{sp}%
%
\special{pn 20}%
\special{pa 3340 360}%
\special{pa 3340 1580}%
\special{fp}%
%
\special{pn 20}%
\special{pa 3660 340}%
\special{pa 3660 1570}%
\special{fp}%
%
\special{pn 20}%
\special{ar 3500 310 150 80  0.0000000 6.2831853}%
%
\special{pn 8}%
\special{ar 3510 4080 560 250  0.0000000 6.2831853}%
%
\special{pn 8}%
\special{pa 4080 330}%
\special{pa 4080 4080}%
\special{fp}%
%
\special{pn 8}%
\special{pa 2950 340}%
\special{pa 2950 4080}%
\special{fp}%
%
\special{pn 8}%
\special{ar 1410 330 560 250  0.0000000 6.2831853}%
%
\special{pn 8}%
\special{ar 1410 4090 560 250  0.0000000 6.2831853}%
%
\special{pn 8}%
\special{pa 1980 340}%
\special{pa 1980 4090}%
\special{fp}%
%
\special{pn 8}%
\special{pa 850 350}%
\special{pa 850 4090}%
\special{fp}%
%
\special{pn 8}%
\special{ar 5890 310 560 250  0.0000000 6.2831853}%
%
\special{pn 8}%
\special{ar 5890 4070 560 250  0.0000000 6.2831853}%
%
\special{pn 8}%
\special{pa 6460 320}%
\special{pa 6460 4070}%
\special{fp}%
%
\special{pn 8}%
\special{pa 5330 330}%
\special{pa 5330 4070}%
\special{fp}%
%
\special{pn 20}%
\special{ar 3500 1590 150 80  0.0000000 6.2831853}%
%
\special{pn 20}%
\special{ar 3500 2710 150 80  0.0000000 6.2831853}%
%
\special{pn 20}%
\special{ar 3520 2070 560 250  0.0000000 6.2831853}%
%
\special{pn 8}%
\special{pa 3660 1600}%
\special{pa 5730 1600}%
\special{dt 0.045}%
\special{pa 5730 1600}%
\special{pa 5729 1600}%
\special{dt 0.045}%
%
\special{pn 8}%
\special{pa 3680 2720}%
\special{pa 5750 2720}%
\special{dt 0.045}%
\special{pa 5750 2720}%
\special{pa 5749 2720}%
\special{dt 0.045}%
%
\special{pn 20}%
\special{ar 5850 1600 150 80  0.0000000 6.2831853}%
%
\special{pn 20}%
\special{ar 5850 2730 150 80  0.0000000 6.2831853}%
%
\special{pn 20}%
\special{pa 5690 1610}%
\special{pa 5690 2700}%
\special{fp}%
%
\special{pn 20}%
\special{pa 6020 1640}%
\special{pa 6020 2730}%
\special{fp}%
\put(12.0000,-46.0000){\makebox(0,0)[lb]{t=-0.5}}%
\put(32.0000,-46.0000){\makebox(0,0)[lb]{t=0}}%
\put(56.0000,-46.0000){\makebox(0,0)[lb]{t=0.5}}%
%
\special{pn 8}%
\special{ar 3510 320 560 250  0.0000000 6.2831853}%
%
\special{pn 8}%
\special{pa 2950 340}%
\special{pa 2950 4080}%
\special{fp}%
%
\special{pn 8}%
\special{pa 4080 330}%
\special{pa 4080 4080}%
\special{fp}%
%
\special{pn 8}%
\special{ar 3510 4080 560 250  0.0000000 6.2831853}%
%
\special{pn 20}%
\special{ar 3500 310 150 80  0.0000000 6.2831853}%
%
\special{pn 20}%
\special{pa 3660 340}%
\special{pa 3660 1570}%
\special{fp}%
%
\special{pn 20}%
\special{pa 3340 360}%
\special{pa 3340 1580}%
\special{fp}%
%
\special{pn 20}%
\special{pa 3350 4000}%
\special{pa 3360 3970}%
\special{pa 3384 3949}%
\special{pa 3413 3935}%
\special{pa 3444 3926}%
\special{pa 3475 3921}%
\special{pa 3507 3920}%
\special{pa 3539 3923}%
\special{pa 3570 3930}%
\special{pa 3601 3940}%
\special{pa 3627 3958}%
\special{pa 3647 3983}%
\special{pa 3648 4014}%
\special{pa 3630 4040}%
\special{pa 3603 4058}%
\special{pa 3573 4070}%
\special{pa 3542 4077}%
\special{pa 3510 4080}%
\special{pa 3478 4079}%
\special{pa 3447 4075}%
\special{pa 3416 4066}%
\special{pa 3387 4053}%
\special{pa 3362 4032}%
\special{pa 3350 4003}%
\special{pa 3350 4000}%
\special{sp}%
%
\special{pn 20}%
\special{pa 3340 3970}%
\special{pa 3340 2740}%
\special{fp}%
%
\special{pn 20}%
\special{pa 3660 3950}%
\special{pa 3660 2730}%
\special{fp}%
%
\special{pn 8}%
\special{ar 1410 330 560 250  0.0000000 6.2831853}%
%
\special{pn 8}%
\special{ar 1410 4090 560 250  0.0000000 6.2831853}%
%
\special{pn 8}%
\special{pa 1980 340}%
\special{pa 1980 4090}%
\special{fp}%
%
\special{pn 8}%
\special{pa 850 350}%
\special{pa 850 4090}%
\special{fp}%
%
\special{pn 8}%
\special{ar 5890 310 560 250  0.0000000 6.2831853}%
%
\special{pn 8}%
\special{ar 5890 4070 560 250  0.0000000 6.2831853}%
%
\special{pn 8}%
\special{pa 6460 320}%
\special{pa 6460 4070}%
\special{fp}%
%
\special{pn 8}%
\special{pa 5330 330}%
\special{pa 5330 4070}%
\special{fp}%
\put(24.3000,-53.3000){\makebox(0,0)[lb]{Figure \ref{ribbonmove}.1}}%
\end{picture}%
\vskip3mm

\unitlength 0.1in
\begin{picture}(56.10,50.60)(8.50,-51.20)
%
\special{pn 8}%
\special{ar 3510 320 560 250  0.0000000 6.2831853}%
%
\special{pn 20}%
\special{pa 3660 3950}%
\special{pa 3660 2730}%
\special{fp}%
%
\special{pn 20}%
\special{pa 3340 3970}%
\special{pa 3340 2740}%
\special{fp}%
%
\special{pn 20}%
\special{pa 3350 4000}%
\special{pa 3360 3970}%
\special{pa 3384 3949}%
\special{pa 3413 3935}%
\special{pa 3444 3926}%
\special{pa 3475 3921}%
\special{pa 3507 3920}%
\special{pa 3539 3923}%
\special{pa 3570 3930}%
\special{pa 3601 3940}%
\special{pa 3627 3958}%
\special{pa 3647 3983}%
\special{pa 3648 4014}%
\special{pa 3630 4040}%
\special{pa 3603 4058}%
\special{pa 3573 4070}%
\special{pa 3542 4077}%
\special{pa 3510 4080}%
\special{pa 3478 4079}%
\special{pa 3447 4075}%
\special{pa 3416 4066}%
\special{pa 3387 4053}%
\special{pa 3362 4032}%
\special{pa 3350 4003}%
\special{pa 3350 4000}%
\special{sp}%
%
\special{pn 20}%
\special{pa 3340 360}%
\special{pa 3340 1580}%
\special{fp}%
%
\special{pn 20}%
\special{pa 3660 340}%
\special{pa 3660 1570}%
\special{fp}%
%
\special{pn 20}%
\special{ar 3500 310 150 80  0.0000000 6.2831853}%
%
\special{pn 8}%
\special{ar 3510 4080 560 250  0.0000000 6.2831853}%
%
\special{pn 8}%
\special{pa 4080 330}%
\special{pa 4080 4080}%
\special{fp}%
%
\special{pn 8}%
\special{pa 2950 340}%
\special{pa 2950 4080}%
\special{fp}%
%
\special{pn 8}%
\special{ar 1410 330 560 250  0.0000000 6.2831853}%
%
\special{pn 8}%
\special{ar 1410 4090 560 250  0.0000000 6.2831853}%
%
\special{pn 8}%
\special{pa 1980 340}%
\special{pa 1980 4090}%
\special{fp}%
%
\special{pn 8}%
\special{pa 850 350}%
\special{pa 850 4090}%
\special{fp}%
%
\special{pn 8}%
\special{ar 5890 310 560 250  0.0000000 6.2831853}%
%
\special{pn 8}%
\special{ar 5890 4070 560 250  0.0000000 6.2831853}%
%
\special{pn 8}%
\special{pa 6460 320}%
\special{pa 6460 4070}%
\special{fp}%
%
\special{pn 8}%
\special{pa 5330 330}%
\special{pa 5330 4070}%
\special{fp}%
%
\special{pn 20}%
\special{ar 3500 1590 150 80  0.0000000 6.2831853}%
%
\special{pn 20}%
\special{ar 3500 2710 150 80  0.0000000 6.2831853}%
%
\special{pn 20}%
\special{ar 3520 2070 560 250  0.0000000 6.2831853}%
\put(12.0000,-46.0000){\makebox(0,0)[lb]{t=-0.5}}%
\put(32.0000,-46.0000){\makebox(0,0)[lb]{t=0}}%
\put(56.0000,-46.0000){\makebox(0,0)[lb]{t=0.5}}%
%
\special{pn 8}%
\special{ar 3510 320 560 250  0.0000000 6.2831853}%
%
\special{pn 8}%
\special{pa 2950 340}%
\special{pa 2950 4080}%
\special{fp}%
%
\special{pn 8}%
\special{pa 4080 330}%
\special{pa 4080 4080}%
\special{fp}%
%
\special{pn 8}%
\special{ar 3510 4080 560 250  0.0000000 6.2831853}%
%
\special{pn 20}%
\special{ar 3500 310 150 80  0.0000000 6.2831853}%
%
\special{pn 20}%
\special{pa 3660 340}%
\special{pa 3660 1570}%
\special{fp}%
%
\special{pn 20}%
\special{pa 3340 360}%
\special{pa 3340 1580}%
\special{fp}%
%
\special{pn 20}%
\special{pa 3350 4000}%
\special{pa 3360 3970}%
\special{pa 3384 3949}%
\special{pa 3413 3935}%
\special{pa 3444 3926}%
\special{pa 3475 3921}%
\special{pa 3507 3920}%
\special{pa 3539 3923}%
\special{pa 3570 3930}%
\special{pa 3601 3940}%
\special{pa 3627 3958}%
\special{pa 3647 3983}%
\special{pa 3648 4014}%
\special{pa 3630 4040}%
\special{pa 3603 4058}%
\special{pa 3573 4070}%
\special{pa 3542 4077}%
\special{pa 3510 4080}%
\special{pa 3478 4079}%
\special{pa 3447 4075}%
\special{pa 3416 4066}%
\special{pa 3387 4053}%
\special{pa 3362 4032}%
\special{pa 3350 4003}%
\special{pa 3350 4000}%
\special{sp}%
%
\special{pn 20}%
\special{pa 3340 3970}%
\special{pa 3340 2740}%
\special{fp}%
%
\special{pn 20}%
\special{pa 3660 3950}%
\special{pa 3660 2730}%
\special{fp}%
%
\special{pn 8}%
\special{ar 1410 330 560 250  0.0000000 6.2831853}%
%
\special{pn 8}%
\special{ar 1410 4090 560 250  0.0000000 6.2831853}%
%
\special{pn 8}%
\special{pa 1980 340}%
\special{pa 1980 4090}%
\special{fp}%
%
\special{pn 8}%
\special{pa 850 350}%
\special{pa 850 4090}%
\special{fp}%
%
\special{pn 8}%
\special{ar 5890 310 560 250  0.0000000 6.2831853}%
%
\special{pn 8}%
\special{ar 5890 4070 560 250  0.0000000 6.2831853}%
%
\special{pn 8}%
\special{pa 6460 320}%
\special{pa 6460 4070}%
\special{fp}%
%
\special{pn 8}%
\special{pa 5330 330}%
\special{pa 5330 4070}%
\special{fp}%
%
\special{pn 8}%
\special{pa 3320 1590}%
\special{pa 1590 1590}%
\special{dt 0.045}%
\special{pa 1590 1590}%
\special{pa 1591 1590}%
\special{dt 0.045}%
%
\special{pn 8}%
\special{pa 3330 2740}%
\special{pa 1620 2740}%
\special{dt 0.045}%
\special{pa 1620 2740}%
\special{pa 1621 2740}%
\special{dt 0.045}%
%
\special{pn 20}%
\special{ar 1420 1590 150 80  0.0000000 6.2831853}%
%
\special{pn 20}%
\special{ar 1420 2720 150 80  0.0000000 6.2831853}%
%
\special{pn 20}%
\special{pa 1260 1600}%
\special{pa 1260 2690}%
\special{fp}%
%
\special{pn 20}%
\special{pa 1590 1630}%
\special{pa 1590 2720}%
\special{fp}%
\put(22.8000,-52.9000){\makebox(0,0)[lb]{Figure \ref{ribbonmove}.2}}%
\end{picture}%
\vskip5mm


We regard $B$ as 
(a close 2-disc)$\times\{s| 0\leqq s\leqq1\}
\times\{t| -1\leqq t\leqq1\}$. 
We let $B_t=$(a close 2-disc)$\times\{s| 0\leqq s\leqq1\}
\times\{t \}$.  
Then $B=\cup B_t$, where $\{t| -1\leqq t\leqq1\}$. 
In Figure \ref{ribbonmove}.1 and \ref{ribbonmove}.2, 
we draw $B_{-0.5}, B_{0}, B_{0.5}$ $\subset B$. 
We draw $K_1$ and $K_2$ by the bold line. 
The fine line denotes $\partial B_t$. 
  
$B\cap K_1$ (resp. $B\cap K_2$) is diffeomorphic to 
$D^2\amalg (S^1\times \{s| 0\leqq s\leqq1\})$, 
where $\amalg$ denotes the disjoint union.

$B\cap K_1$ has the following properties:  
$B_t\cap K_1$ is empty for $-1\leqq t<0$ and $0.5<t\leqq1$.
$B_0\cap K_1$ is diffeomorphic to 
$D^2\amalg(S^1\times \{0\leqq s\leqq0.3\})
\amalg(S^1\times \{0.7\leqq s\leqq1\})$. 
$B_{0.5}\cap K_1$ is diffeomorphic to 
$(S^1\times \{0.3\leqq s\leqq0.7\})$. 
$B_t\cap K_1$ is diffeomorphic to $S^1\amalg S^1$ for $0<t<0.5$. 
(Here we draw $S^1\times \{0\leqq s\leqq1\}$ 
to have the corner in $B_0$ and in $B_{0.5}$. 
Strictly to say, $B\cap K_1$ in $B$ is a smooth embedding 
which is obtained by making the corner smooth naturally.)

$B\cap K_2$ has the following properties:  
$B_t\cap  K_2$ is empty for $-1\leqq t<-0.5$ and $0<t\leqq1$.
$B_0\cap K_2$ is diffeomorphic to 
$D^2\amalg(S^1\times \{0\leqq s\leqq0.3\})
\amalg(S^1\times \{0.7\leqq s\leqq1\})$. 
$B_{-0.5}\cap  K_2$ is diffeomorphic to 
$(S^1\times \{0.3\leqq s\leqq0.7\})$. 
$B_t\cap  K_2$ is diffeomorphic to $S^1\amalg S^1$ for $-0.5<t<0$. 
(Here we draw 
$S^1\times \{0\leqq s\leqq1\}$ to have the corner 
in $B_0$ and in $B_{-0.5}$. 
Strictly to say, $B\cap K_1$ in $B$ is a smooth embedding 
which is obtained by making the corner smooth naturally.)

In Figure \ref{ribbonmove}.1 (resp. \ref{ribbonmove}.2) 
there are an oriented cylinder 
$S^1\times\{0\leqq s\leqq1\}$ 
and an oriented disc $D^2$ as we stated above. 
We do not make any assumption about 
the orientation of the cylinder and the disc. 
The orientation of $B\cap K_1$ (resp. $B\cap K_2$ ) 
coincides with that of the cylinder and that of the disc.



Suppose that $K_2$ is obtained from $K_1$ by one ribbon-move 
and that $K'_2$ is equivalent to $K_2$.   
Then we also say that $K'_2$ is obtained from $K_1$ 
by one {\it ribbon-move}.   
If $K_1$ is obtained from $K_2$ by one ribbon-move,  
then we also say that $K_2$ is obtained from $K_1$ by one {\it ribbon-move}.   
\end{defn}


\vskip3mm
\begin{defn}\label{ribbonmove2}    
Two 2-knots $K_1$ and $K_2$ are said to be {\it ribbon-move equivalent} 
if there are 2-knots 
$K_1=\Hat{K}_1, \Hat{K}_2,...,\Hat{K}_{r-1},\Hat{K}_r=K_2$  
 ($r\in{\N}, p\geq2$) such that 
$\Hat{K}_i$ is obtained from $\Hat{K}_{i-1}$ $(1< i\leqq r)$ 
by one ribbon-move. 
\end{defn}
\vskip3mm

\vskip3mm


\begin{prob}\label{problem}    
Let $K_1$ and $K_2$ be 2-knots.
Find a necessary (resp. sufficient, necessary and sufficient )
condition that $K_1$ and $K_2$ are ribbon-move equivalent.  
\end{prob}
\vskip3mm

In \cite{O1} the author proved the following. 

\vskip3mm
\begin{thm}\label{mu} 
{\rm (\cite{O1})}  
$(1)$ If 2-knots $K$ and $K'$ are ribbon-move equivalent, then 
$$\mu(K)=\mu(K'),$$

\noindent 
where $\mu(\hskip1.5mm)$ denotes the $\mu$-invariant of 2-knots.  

\noindent
$(2)$ Let $K_1$ and $K_2$ be 2-knots in $S^4$. 
Suppose that $K_1$ are ribbon-move equivalent to $K_2$.   
Let $W_i$ be an arbitrary Seifert hypersurface for $K_i$. 
Then the torsion part of 
$\{H_1(W_1)\oplus H_1(W_2)\}$ is congruent to 
$G\oplus G$ for a finite abelian group $G$. 

\noindent
$(3)$ Not all 2-knots are ribbon-move equivalent to the trivial 2-knot.   

\noindent
$(4)$ The converse of $(1)$ is not true.  The converse of $(2)$ is not true. 
\end{thm}
\vskip3mm

\noindent{\bf Note.}  See \cite{O1} for the $\mu$-invariant of 2-knots. 
\vskip3mm

Furthermore, in \cite{O2} the author proved the following Theorem \ref{O2}.

\begin{defn} \label{cano}    
Let $K$ be a 2-knot $\subset S^4$. 
Let $N(K)$ be the tubular neighborhood of $K$ in $S^4$. 
Let $\alpha:\pi_1(\overline{S^4-N(K)})\rightarrow 
H_1(\overline{S^4-N(K)};  {\Z})$ 
be the abelianization. 
Note that any 1-cycle 
is oriented naturally 
by using the orientation of $K$ and that of $S^4$.  
We define the canonical isomorphism  
$\beta:H_1(\overline{S^4-N(K)};  {\Z})\rightarrow  {\Z}$ 
by using this orientation. 
Let $\widetilde{X^\infty_K}$ be 
the covering space associated with 
$\beta\circ\alpha:\pi_1(\overline{S^4-N(K)})\rightarrow  {\Z}$. 
We call $\widetilde{X^\infty_K}$ 
the {\it canonical infinite cyclic covering space} 
of the complement $\overline{S^4-N(K)}$ of $K$.  
We also call $\widetilde{X^\infty_K}$ 
the {\it canonical infinite cyclic covering space} for $K$.
See \cite{Farber78, Levine66, Rolfsen76} 
for canonical infinite cyclic covering spaces for details.
\end{defn}

\begin{note}\label{tub} 
In this paper, 
if we regard a tubular neighborhood as a fiber bundle naturally, 
it is the close disc (not the open disc) that  
the fiber of the fiber bundle  is. 
That is, we have the following. 
Let $A$ be a $a$-submanifold in a $b$-manifold $B$  ($a,b\in\N\cup\{0\}$). 
Then the tubular neighborhood of $A$ in $B$ is 
a fiber bundle over $A$ whose fiber is the $(b-a)$-dimensional close disc.  
\end{note}

\vskip3mm

\begin{thm} $($\cite{O2}$) $\label{O2}    
 Let $K$ and $K'$ be 2-knots. 
Suppose that $K$ and $K'$ are ribbon-move equivalent. 
Then there is an isomorphism  
$$c: {\mathrm {Tor}} H_1(\widetilde{X^\infty_K}; {\Z})\to 
{\mathrm {Tor}} H_1(\widetilde{X^\infty_{K'}}; {\Z}),$$ 
where the homomorphism $c$ is 
not only one as $\Z$-modules but also one as ${\Z}[t,t^{-1}]$-modules, 
with the following properties. 

\noindent$(1)$ 
 Let $x,y\in{\mathrm{Tor}} H_1(\widetilde{X^\infty_K};  {\Z})$. 
Then we have 
$${\mathrm {lk}}(x, y)={\mathrm {lk}}(c(x), c(y)),$$ 
where ${\mathrm {lk}}(\quad)$ denotes the Farber-Levine pairing.  
That is, 
the Farber-Levine pairing on 
${\mathrm{Tor}} H_1(\widetilde{X^\infty_K};  {\Z})$ is equivalent to 
that on 
${\mathrm{Tor}} H_1(\widetilde{X^\infty_{K'}};  {\Z})$. 

\noindent$(2)$ 
 Let $\alpha: H_1(\widetilde{{X}^{\infty}_{K}}; {\Z})\to{\Z_p}$  
be a homomorphism. 
Note that there is a homomorphism 
\newline
$\alpha': H_1(\widetilde{{X}^{\infty}_{K'}}; {\Z})\to{\Z_p}$  
such that 
\newline 
$\alpha\vert_{\mathrm Tor}=(\alpha'\vert_{\mathrm Tor})\circ c$. 
Then we have 
$\tilde\eta(K, \alpha)=\tilde\eta(K', \alpha')\in\Q/\Z$.
That is, 
the set of the values of the $\Q/\Z$-valued $\tilde\eta$ invariants for $K$ 
is equivalent to that for $K'$. 
\end{thm}
\vskip3mm


\noindent{\bf Note.}
See \cite{O2} for the $\widetilde\eta$-invariants of 2-knots,  
the Farber-Levine pairing, and  the Alexander module.

\section{Main results}\label{main}

\begin{defn}\label{oset}
Let $K$ be a 2-knot $\subset S^4$. 
Let $\widetilde{X^\infty_K}$ be the canonical infinite cyclic covering space 
for $K$ (see Definition \ref{cano}). 
Let $\MM=\{M_1,...,M_m\} (m\in\N)$ 
be a set such that 
$M_i$ is an open oriented 3-submanifold $\subset \widetilde{X^\infty_K}$ 
and that 
$[M_i]\in H_3^\infty(\widetilde{X^\infty_K};\R).$
Here, 
$H_*^\infty(\quad;\R)$ 
denotes the $\infty$-chain homology group with the $\R$ coefficient.  
\footnote{
Let $X$ be a topological space. 
The {\it infinity chain homology group} $H^\infty_i(X;\Z)$ 
is defined by using the {\it infinity chain group}  
\nl$C^\infty_i=\{\Sigma_l \nu_l\sigma_l|$ 
$\sigma_l$ is a $i$-simplex. $\nu_l\in\Z$.  
 The number of nonzero $\nu_l$ may be infinite. $\}.$ 

Recall that the homology group $H^\infty_i(X;\Z)$ 
is defined by using the  chain group 
\nl$C_i=\{\Sigma_l \nu_l\sigma_l|$ 
$\sigma_l$ is a $i$-simplex.  $\nu_l\in\Z$. 
 The number of nonzero $\nu_l$ is finite. $\}.$ 

The infinity chain homology group $H^\infty_i(X;\R)$ 
is defined by using the  infinity chain group 
$C^\infty_i\oplus\R$. 

Recall that the homology group $H^\infty_i(X;\R)$
is defined by using the  infinity chain group 
$C_i\oplus\R$. 
}

The set $\MM=\{M_1,...,M_m\}$ is called 
an {\it $o$-set} for $K$ 
if $\MM$ satisfies the following conditions. 

\noindent(1) 
$M_i$ intersects $M_j$ transversely ($i\neq j$). 

\noindent(2) 
$M_i\cap M_j$ intersects $M_k$ transversely 
($i\neq j$, $j\neq k$, and $k\neq i$). 

\noindent(3) 
$M_i\cap M_j\cap M_k$ is 
an oriented submanifold 
$S\amalg R$, where $\amalg$ denotes the disjoint union, 
with the following conditions.

\noindent 
(i) $S$ is a disjoint union of circles. 
We do not assume the number of circles.
The number of circles may be 
$\left\{
\begin{array}{l}
\mbox{zero}\\
\mbox{finite and nonzero}\\
\mbox{infinite.}
\end{array}
\right.$

\noindent
(ii) $R$ is the empty set $\phi$ or 
is diffeomorphic to a single line $\R$ 
with the following property $(*)$. 

$(*)$ 
Let $V$ be any Seifert hypersurface for $K$. 
Let 
$\pi:\widetilde{X^\infty_K}\to\overline{S^4-N(K)}$ be the projection map. 
Then $\pi^{-1}(V)=\amalg_{i=-\infty}^{\infty}V_i$ and 
each $V_i$ is diffeomorphic to $V$. 
Let this $\R$ satisfy the following. 

\noindent (a) 
This $\R$ and $V_i$ intersect transversely. 

\noindent (b) 
The algebraic number of the points, (this $\R$)$\cap V_i$, is one 
for each $i$. 

\noindent (c) 
The geometric number of the points, (this $\R$)$\cap V_i$, is finite  
for each $i$.

\vskip3mm
\noindent 
Note. In this paper, the fact 
 that $\R$ is a submanifold $\subset\widetilde{X^\infty_K}$ 
means the following: 
for each $p\in\R$, there is an open set $U\subset\widetilde{X^\infty_K}$ 
such that $p\in U$ 
and that $\R\cap (U-\R)=\phi$.

Let $\MM=\{M_1,...,M_m\}$ be an $o$-set for $K$. 
Take three elements $M_i, M_j, M_k$. 
Then $M_i\cap M_j\cap M_k$ is $S\amalg R$ as above.
Define  the $o$-invariant 
$$
o(M_i, M_j, M_k)=
\left\{
\begin{array}{ll}
0&\mbox{ if $R$ is the empty set $\phi$}\\
1&\mbox{ if $R\cong\R$.}
\end{array}
\right.
$$
\end{defn}

\vskip3mm

Our main theorem is the following. 
\vskip3mm

\begin{thm}\label{daiji} 
Let $K$ and $K'$ be 2-knots $\subset S^4$. 
Suppose that  $K$ is ribbon-move equivalent to $K'$.  
Let $\MM=\{M_1,...,M_m\} (m\in\N)$ be an $o$-set for $K$. 
Then 
there is an $o$-set 
$\MM'=\{M'_1,...,M'_m\} (m\in\N)$ for $K'$ 
such that 
$o(M_i, M_j, M_k)=o(M'_i, M'_j, M'_k)$ 
for any $i, j, k$. 
\end{thm} 
\vskip3mm

The results in \cite{O1,O2} (Theorem \ref{mu}, \ref{O2} in this paper) 
does not imply Theorem \ref{daiji} (see Note \ref{can}). 
However, this paper prove  Theorem \ref{daiji}.  
Hence Theorem \ref{daiji} is new.

Theorem \ref{daiji} implies Theorem  \ref{cor}.

\vskip3mm
\begin{thm} \label{cor} {\rm (same as Theorem \ref{saisho}.)} 
Let $K$ be a 2-dimensional fibered knot 
whose fiber is the punctured 3-dimensional torus. 
Let $P$ be a 2-dimensional knot whose Seifert hypersurface is 
a punctured homology sphere. 
Then $K$ is not ribbon-move equivalent to $P$. 
\end{thm}
\vskip3mm

\begin{note} \label{Z}
There is a 2-dimensional fibered knot 
whose fiber is the punctured torus (see \cite{AR, AK, CS, Gompf}).
There is a 2-dimensional knot $Z$ whose Seifert hypersurface is 
the punctured Poincar\'e sphere (see \cite{Zeeman}). 
By \cite{O1}, $\mu(Z)\neq0$. 
Take the connected sum $Z\sharp Z$. 
By \cite{O1}, $\mu(Z\sharp Z)=0$. 
Hence there is a 2-knot $J$ whose Seifert hypersurface is a punctured 
homology 3-sphere such that 
$\mu(J)
\left\{
\begin{array}{l}
=0  \\
\neq0. 
\end{array}
\right.
$ 
\end{note}

\vskip3mm
\noindent{\bf Proof of Theorem \ref{cor}.}
By the definition of $K$, 
$\widetilde{X^\infty_K}\cong T^3\times\R$. 
Let $T^3=S^1_1\x S^1_2\x S^1_3.$ 
Let 
$M_1=  p_1\x S^1_2\x S^1_3\times\R$, 
$M_2=S^1_1\x p_2  \x S^1_3\times\R$, and 
$M_3=S^1_1\x S^1_2\x p_3  \times\R$, 
where  
$p_i$ is a point $\in S^1_i$. 
Then we have the following. 

\noindent
(1) 
$\{M_1, M_2, M_3\}$ is an $o$-set. 

\noindent
(2) 
$M_1\cap M_2\cap M_3=\R$. 

\noindent
(3) 
$o(M_1, M_2, M_3)=1$.

 We prove by contradiction. 
 We suppose that 
$K$ is ribbon-move equivalent to $P$. 
By Theorem \ref{daiji}, 
there is an  $o$-set 
$\{M'_1, M'_2, M'_3\}$ for the 2-knot $P$ such that 
$o(M'_1, M'_2, M'_3)=1$. 
Hence $M'_1\cap M'_2\cap M'_3$ represents a nonzero element 
$\in H_1^\infty(\widetilde{X^\infty_P};\R)$.
Hence 
$[M'_i] (i=1,2,3)$ represents a nonzero element 
$\in H_3^\infty(\widetilde{X^\infty_P};\R)$.
Hence 
$H_3^\infty(\widetilde{X^\infty_P};\R)$ is not congruent to 0. 
However, since a punctured homology sphere is a Seifert hypersurface for 
the 2-knot $P$, $H_3^\infty(\widetilde{X^\infty_P};\R)\cong0$. 
We arrived at a contradiction. Hence the initial condition is false. 
That is,  $K$ is not ribbon-move equivalent to $P$. 

Note that, only in the $\mu(K)\neq\mu(P)$ case, 
Theorem \ref{mu} implies that $K$ is not ribbon equivalent to $P$. 

\vskip3mm

\begin{note}\label{trivial}
In the similar manner in the above proof, 
we have the following. 
Suppose that $\{M_1, M_2, M_3\}$ is an $o$-set for the trivial 2-knot. 
Then $o(M_1, M_2, M_3)=0$. 
\end{note}

\vskip3mm

\begin{note}\label{can}
Theorem \ref{mu}, \ref{O2} cannot prove 
the $\mu(K)=\mu(P)$ case of Theorem \ref{cor}. 
However, Theorem \ref{daiji} (and Theorem \ref{cor}) can.
\end{note}

\vskip3mm

We have the following Theorem \ref{gyaku}. 
Compare Theorem \ref{gyaku} to Note \ref{trivial}. 
Hence the converse of Theorem \ref{daiji} 
is not true.

\vskip3mm
\begin{thm}\label{gyaku} 
There is a 2-knot $K$ with the following properties.

\noindent(1) 
For any 
$o$-set $\{M_1, M_2, M_3\}$ for $K$, 
$o(M_1, M_2, M_3)=0$. 
 
\noindent(2) 
$K$ is not ribbon-move equivalent to the trivial knot. 
\end{thm}

\noindent{\bf Proof.} 
The 2-knot $Z$ in Note \ref{Z} is an example because of Theorem 2.4.(1).

\vskip3mm

Compare the following theorem 
to the above Theorem \ref{daiji} (resp. \ref{cor}).  
\vskip3mm

\begin{thm}\label{ichini}
There is a nontrivial 2-knot K with the following properties.

\noindent(1) 
$K$ is ribbon-move equivalent to the trivial 2-knot.

\noindent(2) 
There are 
a 3-dimensional open oriented submanifold $M$ 
such that 
$[M]\in H_3^\infty(\widetilde{X^\infty_K};\R)$ 
and 
a 2-dimensional open oriented submanifold $N$ 
such that 
$[N]\in H_2^\infty(\widetilde{X^\infty_K};\R)$ 
with the following properties.

\noindent$(i)$  
$M\cap N\cong \R$.

\noindent$(ii)$  
This $\R$ satisfies the $(*)$ in (3)(ii) in Definition \ref{oset}. 
\end{thm}

\vskip3mm
\noindent
{\bf Proof of Theorem \ref{ichini}.} 
An example is the spun knot of the trefoil knot.

\section{(1,2)-pass-moves and ribbon-move surgeries of $S^4$}\label{12}

\noindent 
In order to prove our main theorem (Theorem \ref{daiji}), 
we use the (1,2)-pass-moves for 2-knots. 
\cite{O1} defined the (1,2)-pass-moves for 2-knots.

\vskip3mm


\begin{defn}\label{passmove}     
Let $K_1$ and $K_2$ be 2-links in $S^4$. 
We say that $K_2$ is obtained from $K_1$ by one {\it (1,2)-pass-move } 
if there is a 4-ball $B$ embedded in $S^4$ with the following properties.  

We draw $B$ as in Definition 1.1. 

\vskip1mm
\noindent
(1) $K_1$ coincides with $K_2$ in 
$\overline{S^4-B}$. 
This identity map from 
$\overline{K_1-B}$ to $\overline{K_2-B}$  
is orientation preserving. 
Note that this condition on $K_i$ implies that 
$K_1$  coincides with $K_2$ in $\overline{S^4-B}$.

\vskip1mm\noindent
(2) $B\cap K_1$  is drawn as in 
Figure \ref{passmove}.1. 
$B\cap K_2$ is drawn as in 
Figure \ref{passmove}.2.  

We suppose that 
each vector $\overrightarrow{x}$, $\overrightarrow{y}$ 
in Figure \ref{passmove}.1 (resp. \ref{passmove}.2)   
is a tangent vector of each disc at a point. 
(Note we use $\overrightarrow{x}$ (resp. $\overrightarrow{y}$)
for different vectors.)
The orientation of each disc 
in 
Figure \ref{passmove}.1 (resp. Figure \ref{passmove}.2)   
is determined by the each set 
$\{\overrightarrow{x},\overrightarrow{y}\}$. 
We do not make any assumption about 
the orientations of the annuli 
in the 
Figure \ref{passmove}.1 (resp. Figure \ref{passmove}.2).   
The orientation of $B\cap K_1$ coincides with that of 
the disjoint union of the two discs and the annuli. 

\np
\vskip3mm

\input 6.1.tex

\vskip3mm

\input 6.2.tex
\vskip3mm

Suppose that $K_2$ is obtained from $K_1$ by one (1,2)-pass-move 
and that $K'_2$ is equivalent to $K_2$.   
Then we also say that $K'_2$ is obtained from $K_1$ by one 
{\it (1,2)-pass-move }.

If $K_1$ is obtained from $K_2$ by one (1,2)-pass-move,  then 
we also say that $K_2$ is obtained from $K_1$ by one {\it (1,2)-pass-move }.

2-links $K_1$ and $K_2$ are said to be {\it (1,2)-pass-move equivalent} 
if there are 2-links 
$K_1=\Hat{K}_1, \Hat{K}_2,...,\Hat{K}_{p-1},\Hat{K}_p=L_2$  
$( p\in\mathbb{N}, p\geq2 )$ such that 
$\Hat{K}_i$ is obtained from $\Hat{K}_{i-1}$ $(1< i\leqq p)$ 
by one (1,2)-pass-move. 
\end{defn}

\vskip3mm
\noindent{\bf Note.} 
In \cite{O10} the author defined ($p,q$)-pass-moves for $p,q\in\N$. 
The (1,2)-pass move here is the ($p,q$)-pass-move in the 
$p=1$ and $q=2$ case there.  

Before \cite{O10}, the author defined other local moves in \cite{O4}. 
The local moves in \cite{O4} are the $(p,p)$-pass-moves in \cite{O10}. 
\vskip3mm

\vskip3mm 
\cite{O1} proved: 

\vskip3mm

\begin{thm}\label{ribbonpass} 
{\rm (\cite{O1})}   
Let $K$ and $K'$ be 2-knots.  
The following conditions $(1)$ and $(2)$ are equivalent. 

$(1)$ $K$ is  (1,2)-pass-move equivalent to $K'$.

$(2)$ 
$K$ is  ribbon-move equivalent to $K'$. 

Furthermore, 
if $K$ is obtained from $K'$ by one ribbon-move, 
then $K$ is obtained from $K'$ by one (1,2)-pass-move.
\end{thm}
\vskip3mm

Let $K_<$ and $K_>$ be 2-knots. 
Suppose that $K_<$ is ribbon-move equivalent to $K_>$. 
By Theorem \ref{ribbonpass}, 
$K_<$ is (1,2)-pass-move equivalent to $K_>$.  
Therefore, 
in order to prove our main result (Theorem \ref{daiji}) in \S\ref{main},  
it suffices to prove 
the case when 
$K'$ is obtained from $K$ 
by one (1,2)-pass-move in a 4-ball $B$ embedded in $S^4$.

\vskip3mm
Next we state a relation among  surgeries, (1,2)-pass-moves and ribbon-moves.

\begin{defn}\label{sur}
Let $M$ be an $m$-manifold or an $m$-manifold with boundary. 
Make a product manifold $M\x [0,1]$. 
Identify $M\x\{0\}$ with $M$. 
Suppose that 
 handles $h^p$ 
 are attached to $M\x [0,1]$, 
where $h^p\cap(M\x [0,1])\subset M\x\{1\}$. 
\newline Let 
$M'=(\partial(h^p\cup(M\x[0,1])))-M\x\{0\}-(\partial M)\x[0,1)$, 
where $(\partial M)\x[0,1)=\phi$ if $\partial M=\phi$. 

If we do the above procedure, 
we say that 
{\it $M'$ is obtained from $M$ 
by the surgery by using the above handles}. 
(In other words, 
if we do the above procedure, 
we sometimes do not explain that we use $M\x[0,1]$.)

See \cite{Smale} for surgeries. 
\end{defn}

\begin{thm}\label{PQ}  
Let $K$ be a 2-knot $\subset S^4$. 
Take a 4-ball $B\subset S^4$.  
Let $K\cap B$ be as in Figure \ref{passmove}.1.
Then there is a submanifold 
$P\amalg Q\subset B$ with the following properties.  

\noindent(1) 
$P\cong S^1$. 

\noindent(2) 
$Q\cong S^2$. 

\noindent(3) 
$P\cap Q=\phi$.  
The linking number of $P$ and $Q$ is one 
if we give an orientation to $P\amalg Q$. 

\noindent(4) 
Regard 
$\left\{
\begin{array}{l}
B\\
S^4
\end{array}
\right.
$ 
as 
$
\left\{
\begin{array}{l}
B\x\{1\}\\
S^4\x\{1\}
\end{array}
\right.
$ 
$\subset
\left\{
\begin{array}{l}
B\x[0,1]\\ 
S^4\x[0,1].  
\end{array} 
\right. 
$

Attach a 5-dimensional 2-handle along $P$ 
with the trivial framing 
and 
a 5-dimensional 3-handle along $Q$ 
with the trivial framing. 
Then the new 4-manifold made from 
$
\left\{
\begin{array}{l}
B\x\{1\}\\
S^4\x\{1\}
\end{array}
\right.
$ 
by this surgery 
is diffeomorphic to 
$\left\{
\begin{array}{l}
B\\
S^4
\end{array}
\right.
$ 
again.

\noindent(5) 
The new knot in the new $S^4$ obtained by these surgeries 
is the knot made from $K$ by one (1,2)-pass-move in $B$.  
That is,  $K\cap$ (the new $B$) is as in Figure \ref{passmove}.2.

\end{thm}

\vskip3mm 
\noindent{\bf Note.} 
By Theorem \ref{ribbonpass}, the new knot is ribbon-move equivalent to $K$. 
\vskip3mm

\begin{defn}\label{surgeryribbon} 
The set of these surgeries in Theorem \ref{PQ} 
is called 
the {\it ribbon move surgery} of $S^4$ along $P\amalg Q$.  
\end{defn}

\np
\noindent{\bf Proof of Theorem \ref{PQ}.} 
See Figure \ref{surgeryribbon}.1.

\vskip3mm
\input figuresurgery.tex   
\vskip5mm

Note that $K\cap B$ is a disjoint union of the cylinder $A$ and 
a set $D$ of the two discs.  
Take a 3-dimensional 1-handle $h^1$ trivially embedded in $B$. 
The manifold $\overline{\partial h^1-\partial B}$ can be regarded as $A$. 
Take a 4-dimensional 2-handle $h^2$ trivially embedded in $B$. 
The manifold $\overline{\partial h^2-\partial B}$ can be regarded as $D$. 

Embed $S^1$ in $B$ 
so that the linking number  of $S^1$ and $h^2$ in $B$  is one, 
let $P$ denote this $S^1$. 
Embed $S^2$ in $B$ 
so that the linking number  of $S^2$ and $h^1$ in $B$  is one, 
let $Q$ denote this $S^2$. 
Note that 
we can define 
the linking number of 
$\left\{
\begin{array}{l} 
P\amalg h^2 \\
Q\amalg h^1 
\end{array}
\right.
$   
in $B$ as above 
if we give an orientation to 
$\left\{
\begin{array}{l} 
P\amalg h^2 \\
Q\amalg h^1. 
\end{array}
\right.
$   
Note that the attaching part of $h^*$ is fixed at $\partial B$. 

Then we can suppose that  the linking number  
of $P (\cong S^1)$ and $Q (\cong S^2)$ in $B$  is one
if we give an orientation to 
$P\amalg h^2\amalg Q\amalg h^1$
so that we do not change 
the linking number of 
$\left\{
\begin{array}{l} 
P\amalg h^2 \\
Q\amalg h^1 
\end{array}
\right.
$   
in $B$. 
This completes the proof.

\vskip1cm
\begin{note}\label{Pont}  
Take a Seifert hypersurface $V$ for $K$. 
We can suppose that $V\cap B=h^1\amalg h^2$. 
Proof. The idea of the proof is Pontrjagin-Thom construction. 
See e.g. P. 49  of \cite{Kirby89}.    
\end{note}

\section{ Ribbon-move surgeries of 
canonical infinite cyclic covering spaces
 $\widetilde{X^\infty_K}$}\label{ccc}

Let $X_K$ denote $\overline{S^4-N(K)}$. 
Note that $\overline{N(K)}$ is diffeomorphic to $K\x$ (a close 2-disc). 
Let $V$ be a Seifert hypersurface for $K$. 
Let the manifold $V\cap X_K$ be called $V$ again. 
Let $N(V)$ be the tubular neighborhood of $V$ in $X_K$. 
Then $\overline{N(V)}\cong V\x[-1,1]$. 
Let $Y_K=\overline{X_K-N(V)}$. 

Let $NV$ denote $\overline{N(V)}$. 
Let $\widetilde{X^\infty_K}$ be the canonical infinite cyclic covering space 
for $K$ (see Definition \ref{cano}). 
We can regard 
$\widetilde{X^\infty_K}=\cup^\infty_{i=-\infty}(NV_i\cup Y_i)$ ($i\in\Z$), 
where 
$\left\{
\begin{array}{l} 
NV_i \\
Y_{Ki} 
\end{array}
\right.
$ 
is made from 
$\left\{
\begin{array}{l}
NV=V\x[-1,1] \\
Y_{K}
\end{array}
\right.$ 
naturally. 
We sometimes abbreviate 
$\left\{
\begin{array}{l}
Y_{Ki} \\
Y_{K}
\end{array}
\right.$ 
to 
$\left\{
\begin{array}{l}
Y_{i} \\
Y. 
\end{array}
\right.$ 
We can regard $NV_i=V_i\x[-1,1].$

\noindent 
Suppose that 
$Y_j\cap NV_i$
$\left\{
\begin{array}{lll}
\neq\phi& {\rm if}&   i=j,  i=j+1 \\
=\phi   & {\rm if}&  i\neq j, i\neq j+1. 
\end{array}
\right.$

Note we can make 
$X_K$ from $NV_i$ and 
$\left\{
\begin{array}{l}
Y_i \\
Y_{i+1}
\end{array}
\right.$ 
by using an attaching map 
which is different from the attaching map in making $X^\infty_K$.

Take $P\amalg Q\subset B$ as in Definition \ref{surgeryribbon}.
We define 
$\left\{
\begin{array}{l}
P_i\\
Q_i
\end{array}
\right.$ 
so that 
$\left\{
\begin{array}{l}
\amalg_{i=-\infty}^{\infty}P_i\\
\amalg_{i=-\infty}^{\infty}Q_i
\end{array}
\right.$ 
is the lift of 
$\left\{
\begin{array}{l}
P\\
Q. 
\end{array}
\right.$ 

\noindent
We define that 
$\left\{
\begin{array}{l}
(\amalg^{\infty}_{j=-\infty}P_j)\cap Y_i\\
(\amalg^{\infty}_{j=-\infty}Q_j)\cap Y_i
\end{array}
\right.$ 
is 
$\left\{
\begin{array}{l}
P_i\\
Q_i. 
\end{array}
\right.$

Note that we can let $P_i\cap NV_j=\phi$ and $Q_i\cap NV_j=\phi$ 
for each $(i,j)$.

\vskip3mm

We use the following kind of surgeries from now on.

\begin{defn}\label{subsur}
Let $N$ be an $n$-submanifold or 
an `$n$-submanifold with boundary' 
of an $m$-manifold $M$. 
Let $N\x [0,1]$ be a submanifold $\subset M$. 
Identify $N\x\{0\}$ with $N$. 
Suppose that 
 handles $h^p$ which are embedded in $M$ 
are attached to $N\x [0,1]$, 
where $h^p\cap(N\x [0,1])\subset N\x\{1\}$. 
Suppose that $h^p\cap (N\x[0,1])=\phi$.  
\newline Let 
$N'=(\partial(h^p\cup(N\x[0,1])))-N\x\{0\}-(\partial N)\x[0,1)$, 
where $(\partial N)\x[0,1)=\phi$ if $\partial N=\phi$. 

If we do the above procedure, 
we say that 
{\it $N'$ is obtained from $N$ 
by the surgery by using the above embedded handles}. 
(In other words, 
if we do the above procedure, 
we sometimes do not explain that we use $N\x[0,1]$.)
\end{defn}

\vskip3mm

\noindent
{\bf Note.} 
We sometimes abbreviate `submanifold with boundary' to submanifold. 


\vskip3mm
Compare Definition \ref{subsur} to Definition \ref{sur}.

\vskip3mm 
\noindent 
\begin{cla}\label{coveringsurgery}
There are surgeries along $P_i\amalg Q_i(\subset Y_{Ki})$ 
with the following property (*).

\noindent(*) 
Let $Y'_{Ki}$ denote the new manifold made from $Y_{Ki}$. 
Let $K'$ be obtained from $K$ 
by the one ribbon-move in $B$. 
It holds that 
$Y'_{Ki}$ is $Y_{K'i}$.  
It holds that 
$\amalg^{\infty}_{i=-\infty}(NV_i\cup Y'_{Ki})
=\widetilde{X^\infty_{K'}}$. 
\end{cla}
\vskip3mm

\vskip3mm
\noindent
{\bf Note.} 
$
\left\{
\begin{array}{l}
P_i \\
Q_i
\end{array}
\right.
$ 
may be a nonvanishing cycle. 

\noindent
\begin{defn}\label{surgeryribboncovering}
The set of these surgeries is called the {\it ribbon-move surgery}  
of the infinite cyclic covering space $\widetilde{X^\infty_{K}}$ 
along $P_i\amalg Q_i$. 
\end{defn}

Compare Definition \ref{surgeryribboncovering} to Definition \ref{surgeryribbon} 
\vskip3mm

\noindent
{\bf Proof of Claim \ref{coveringsurgery}.} 
Carry out the ribbon surgery along $P\amalg Q (\subset Y_K \subset S^4)$. 
Corresponding to this ribbon-surgery, 
we can carry out surgeries on 
$Y_{Ki}$ along $P_i\amalg Q_i$ 
to satisfy the conditions in Claim \ref{coveringsurgery}.

\section{Proof of Theorem \ref{daiji} (main theorem)}\label{proof}


\noindent{\bf Proof of Theorem \ref{daiji}.} 
As we state right after Theorem \ref{ribbonpass}, 
it suffices to prove  
 the case where a 2-knot $K$ is obtained from a 2-knot $K'$ 
by using a single (1,2)-pass-move. 

Let 
$\MM=\{\widetilde{M_1},...,\widetilde{M_m}\} (m\in\N)$ 
be an $o$-set for $K$. 

We can suppose that $\widetilde{M_*}$ and 
$
\left\{
\begin{array}{l}
NV_i \\
V_i
\end{array}
\right.
$ 
intersect transversely.

Note that $\widetilde{M_*}$ represents an element 
$\in H^\infty_3(\widetilde{X^\infty_K};\R)$. 
Recall that 
$V_i$ is a compact oriented manifold 
and 
that $\partial V_i\subset \partial(\widetilde{X^\infty_K})$.

Hence 
$
\left\{
\begin{array}{l}
\widetilde{M_*}\cap NV_i \\
\widetilde{M_*}\cap V_i 
\end{array} 
\right.
$ 
is a compact oriented manifold 
$
\left\{
\begin{array}{l}
F\x[-1,1] \\
F, 
\end{array}
\right.
$ 
where the following hold. 

\noindent (1) 
$F$ is a closed oriented surface. ($F$ may not be connected.)

\noindent (2) 
Let $e_t$ be the embedding map 
$F\x\{t\}\hookrightarrow V_i\x\{t\}$. 
Then $e_t$ and $e_{t'}$ are the same embedding maps 
(
$
\left\{
\begin{array}{l}
-1\leqq t\leqq1 \\
-1\leqq t'\leqq1
\end{array}
\right.
$ 
) 
if we identify 
$F\x\{t\}$ with $F\x\{t'\}$ 
and 
if we do 
$V_i\x\{t\}$ with $V_i\x\{t'\}$.

\noindent 
Let $\widetilde{M_*}\cap NV_i=F\x[-1,1]$  be called $M_{*NVi}$.

Take $h^1\subset V_i$ 
as in Proof of Theorem \ref{PQ}. 
Consider $h^1\x[-1,1]\subset NV_i$.

\vskip3mm

\noindent  
\begin{cla}\label{h1saketa} 
We can suppose that 
the manifold 
$F\times[-1,1]=\widetilde{M_*}\cap NV_i$
is embedded in 
the manifold 
$NV_i-(h^1\x[-1,1])$.  
\end{cla}
\vskip3mm

\noindent  Proof. 
Take the cocore $C$ of $h^1\subset V_i$. 
Note $C$ is a 2-disc. 
We can suppose that 
$F$ and $C$ intersect transversely, if they intersect. 
The intersection $F\cap C$ is a set of circles. 
These circles are the boundaries of discs $\Hat{D}$
such that 
these discs $\Hat{D}$ are embedded in $C$.  
(Note that the discs $\Hat{D}$ 
may intersect each other but our proof may not mind that.)
Regard these 2-discs $\Hat{D}$ as the cores of 4-dimnsional 2-handles. 
Carrying out surgeries on $\widetilde{M_*}$  
by attaching these 2-handles $\Hat{D}$ 
along $F\cap C$, 
let $\widetilde{M_*}\cap (h^1\x[0,1])=\phi$ 
and 
 $(\widetilde{M_*}\cap V_i)\cap h^1=\phi$.  

\noindent 
Note: If the discs $\Hat{D}$ intersect each other, 
use the isotopy of $\Hat{D}$.

Here, 
if necessary,
we carry out these surgeries on $\widetilde{M_*} (*=1,...,m)$. 
These surgeries are  done 
in the interior of 
the tubular neighborhood \nl$N(h^1\x[-1,1])$ of $h^1\x[-1,1]$ 
in $\widetilde{X^\infty_K}$.
Note that 
\nl$N(h^1\x[-1,1])$ is a compact set.  (See Note \ref{tub}.)
Suppose that 
\nl$o(\widetilde{M_\alpha}, \widetilde{M_\beta}, \widetilde{M_\gamma})=1$ 
for a set $\{\alpha, \beta, \gamma\}$. 
By using isotopy of $\widetilde{M_*}$, 
we  can suppose that 
$\widetilde{M_\alpha}\cap \widetilde{M_\beta}\cap \widetilde{M_\gamma}$
does not intersect $N(h^1\x[-1,1])$ before these surgeries. 

Then we have the following. 
Suppose that we obtain  a new intersection 
$\widetilde{M_\alpha}\cap \widetilde{M_\beta}\cap \widetilde{M_\gamma}$ 
after these surgeries. 
(Note that the new intersection consists of triple points 
and is an oriented 1-dimensional manifold.) 
Then the new intersection is a disjoint union of circles
because of the following.

\noindent 
(1) 
These surgeries are done in the interior of the compact set $N(h^1\x[-1,1])$. 

\noindent 
(2) 
The triple point set 
does not exist 
in $N(h^1\x[-1,1])$ 
before these surgeries. 
(Because of the above procedure.)

\vskip3mm
Therefore, even if we change $\widetilde{M_*}$ in the above procedure,   
the $o$-invariant $o(\widetilde{M_\alpha}, \widetilde{M_\beta}, \widetilde{M_\gamma})$ 
for an arbitrary set $\{\alpha, \beta, \gamma\}$ does not change. 

This completes the proof. 


\vskip3mm

Suppose that 
$o(\widetilde{M_\alpha}, \widetilde{M_\beta}, \widetilde{M_\gamma})=1$ 
for a set $\{\alpha, \beta, \gamma\}$. 
Then 
\nl$\widetilde{M_\alpha}\cap \widetilde{M_\beta}\cap \widetilde{M_\gamma}=\R$.

\begin{cla}\label{geometricone}
The geometric number of the points, (this $\R$)$\cap V_i$,  is one 
for each $i$. 
\end{cla}

\noindent 
Proof. 
Recall Definition \ref{oset}.
The algebraic number of the points, (this $\R$)$\cap V_i$,  is one 
for each $i$. 
The geometric number of the points, (this $\R$)$\cap V_i$,  is finite 
for each $i$.

We can suppose that (this $\R$)$\cap P_i=\phi$ for each $i$. 
Because: If it is not an empty set, use the isotopy of $P_i$.

We can suppose that (this $\R$)$\cap Q_i=\phi$ for each $i$. 
Because: If it is not an empty set, use the isotopy of $Q_i$.

If the geometric number, (this $\R$)$\cap V_i$,  is not one, 
carry out surgeries on $V_i$ by using 4-dimensional 1-handles 
with the following properties.

\noindent 
(1) Each of the handles is embedded in 
$\widetilde{X^\infty_K}$.

\noindent
(2) The core of each of the handles is 
a 1-dimensional `submanifold with boundary' of 
this $\R$.

Note that, 
since the geometric number of the points, (this $\R$)$\cap V_i$,  is finite, 
the number of these surgeries is finite.

Note that
the new $V_i$ is orientable.  Because: 
We can suppose that 
the two points 
along which each of the above 1-handles is attached 
have the opposite orientations. 

\vskip3mm
This completes the proof. 

\vskip3mm
We can suppose that 
$\widetilde{M_*}$  and $Y_i$ intersect transversely. 

Let $M_{*i}=Y_i\cap \widetilde{M_*}$ ($*=1,..,m$). 
Then we have the following. 

\noindent(1) 
$\partial M_{*i}=N_{*i-}\amalg N_{*i+}\subset \partial Y_i$. 

\noindent(2) 
$M_{*i}\cap NV_i=N_{*i-}$. 

\noindent(3) 
$M_{*i}\cap NV_{i+1}=N_{*i+}$. 

\noindent(4) 
For each set $\{\alpha, \beta, \gamma\}\subset \{1,...,m\}$, 
$M_{\alpha i}\cap M_{\beta i}\cap M_{\gamma i}$ is 
a disjoint union $S\amalg I$, 
where $S$ is a set of circles 
and where $I$ is the empty set or a single segment $I$ 
with the following properties: 

\noindent
(i) (One of $\partial I)$$\subset NV_i$. 

\noindent
(ii) (The other of $\partial I)$$\subset NV_{i+1}$. 

\noindent
Note that $\partial I$ is two points.
Note that   
($\widetilde{M_\alpha}\cap \widetilde{M_\beta}\cap \widetilde{M_\gamma})
\cap V_i$ 
is geometrically one point. 
(See Claim \ref{geometricone}.)

\begin{cla}
We can suppose that $Q_i\cap M_{*i}=\phi$. 
\end{cla}


\noindent 
{\bf Proof.} 
See Figure \ref{surgeryribbon}.1. 
We can take $Q$ in the tubular neighborhood of $h^1$ in $S^4$. 
Let $C(\widetilde{X^\infty_K})$ be the collar neighborhood of
$\widetilde{X^\infty_K}$. 
Take the above  $N(h^1\x[-1,1])$, 
which is the tubular neighborhood of 
$h^1\x[-1,1]$ in $\widetilde{X^\infty_K}$. 

Then it holds that we can take $Q_i$ 
in $U=N(h^1\x[-1,1])\cup C(\widetilde{X^\infty_K})$.

By the definition of $\widetilde{M_i}$ and Claim \ref{h1saketa},  
we can take $\widetilde{M_i}$ outside $U$. 
Hence $Q_i\cap M_{*i}=\phi$.


\vskip3mm
Next we consider 
$P_i\cap M_{*i}$. (Recall that 
$P_i\cong S^1$.)
We can suppose that 
$P_i\cap M_{*i}$ is a finite set of points. 
Each point is oriented by 
$P_i, M_{*i},$ and $S^4$. 

\vskip3mm  
\noindent Claim. {\it We can suppose that the orientations of these points are same. }
\vskip3mm

\noindent
Proof. Suppose that there are two points 
($\subset P_i\cap M_{*i}$)
such that the two points sits side by side in $P_i$ 
and that they have different orientations.

Carry out a surgery on $M_{*i}$ 
along the two points 
by using a 4-dimensioanl 1-handle with the following properties. 

\noindent 
(1) The handle is embedded in 
$Y_i$. 

\noindent
(2) The core of the handle is 
a 1-dimensional `submanifold with boundary' of 
$P_i$.  
(Note $P_i\cong S^1$.)

\vskip3mm
Note that the new $M_{*i}$ is orientable. 
Because: 
The two points 
along which the above 1-handle is attached 
have the opposite orientations.

Repeat this surgery. 

If necessary, we carry out these surgeries on all $M_{*i}(i=1,...,m)$.  
These surgeries are done 
in the interior of the compact set 
$$W=\overline{N(P_i)\cap \{N(h^1\x[-1,1])-C(\widetilde{X^\infty_K})\}},$$ 
where we have the following.

\noindent(1) $N(P_i)$ is the tubular neighborhood of $P_i$ in 
$Y_i$. 

\noindent(2) We take $N(h^1\x[-1,1])$ so that 
\nl$P_i\cap$(the interior of $N(h^1\x[-1,1])$)$\neq\phi$ 
and that 
\nl$P_i\cap$(the interior of $N(h^1\x[-1,1])$) is connected.

\vskip3mm
\noindent
Suppose that  the triple point set 
$\widetilde{M_\alpha}\cap \widetilde{M_\beta}\cap \widetilde{M_\gamma}$ 
intersects $W$ before these surgeries. 
By using isotopy of $\widetilde{M_*}$, 
we  can suppose that 
$\widetilde{M_\alpha}\cap \widetilde{M_\beta}\cap \widetilde{M_\gamma}$
does not intersect $W$. 
Suppose that 
the triple point set 
$\widetilde{M_\alpha}\cap \widetilde{M_\beta}\cap \widetilde{M_\gamma}$ 
is obtained after these surgeries. 
Then 
the new triple point set 
$\widetilde{M_\alpha}\cap \widetilde{M_\beta}\cap \widetilde{M_\gamma}$ 
is a disjoint union of circles. 
Because we have the following. 

\noindent 
(1) 
These surgeries are done in the interior of the compact set 
$W$.

\noindent 
(2)
The triple point set does not exist in $W$ before these surgeries. 
(Because of the above procedure.)  

\vskip3mm

Therefore, 
the $o$-invariant $o(\widetilde{M_\alpha}, \widetilde{M_\beta}, \widetilde{M_\gamma})$ 
for an arbitrary set $\{\alpha, \beta, \gamma\}$ does not change. 
This completes the proof.

\vskip3mm

\begin{cla}\label{E} 
We can suppose that 
there is a compact 3-dimensional 
`submanifold with boundary'  
$E\subset X_K = \overline{S^4-N(K)}$ 
such that 

\noindent 
$E\cong ($ the punctured $S^1\x S^2)$   
and that 
$\partial E=Q$. 
(Recall $Q\cong S^2$.)
\end{cla} 

\vskip3mm

\noindent{\bf Proof.}  
We can take $Q$ in $t=0 $ in Figure \ref{passmove}.1. 
Take a 3-ball $G$ in $t=0 $ in Figure \ref{passmove}.1 so that $\partial G=Q$.
Note that 
$G\cap K$ is a single circle $S^1$. 


Attach a 4-dimensional 2-handle $h^2$ to $G$ along this $S^1$, 
where the following hold.

\noindent 
(1) $h^2$ is embedded in $S^4$.

\noindent
(2)the core of $h^2$ is a 2-dimensional 
`submanifold with boundary' of $K$.

\vskip3mm
After this surgery by this $h^2$, 
$G$ is changed into $E$ as in Claim \ref{E}. 
This completes the proof. 

\vskip3mm
Note that $E\cap P$ is a single point by the construction 
(see Figure 4.1.1 and Figure \ref{surgeryribbon}.1). 
Note that $E\cap V\neq\phi$.

\vskip3mm
Let $\widetilde{E}$ be the lift of $E$ 
associated with the projection map 
$\widetilde{X^\infty_K}\to X_K$. 

Let $E_{Yi}=\widetilde{E}\cap Y_i$. 
Let $E_{NVi}=\widetilde{E}\cap NV_i$.

Note that $\partial E_{Yi}-\partial Y_i=Q_i$ 
by the construction 
(see Figure 4.1.1 and Figure \ref{surgeryribbon}.1). 
Note that we do not suppose that 
the orientation of $\partial E_{Yi}$ 
coincides with that of $Q_i$. Their orientations may coincide or may not. 
We determine the orientation of 
$E_{Yi}$ (and hence that of $\partial E_{Yi}-\partial Y_i$) 
after several lines from here. 
Take the tubular neighborhood $N(Q_i)$ of $Q_i$ in 
$Y_i$. 
Let us $\overline{E_{Yi}-N(Q_i)}$ call $E_{Yi}$ again.  

\vskip3mm
\noindent 
Claim. {\it We have $E_{Yi}\cap M_{*i}=\phi.$}
\vskip3mm

\noindent  
Proof. 
We can suppose that 
$E_{Yi}\subset N(h^1\x[-1,1])\cap C(\widetilde{X^\infty_K})$. 
We can suppose that 
$M_{*i}$ exists outside 
$N(h^1\x[-1,1])\cap C(\widetilde{X^\infty_K})$.
Therefore  $E_{Yi}\cap M_{*i}=\phi.$

\vskip3mm

Let $\nu$ be the number of the points 
$P_i\cap M_{*i}$ ($\nu\in\{0\}\cup\N$). 
Take $\nu$ copies of $E_{Yi}$. 
Let each $E_{Yi}$ be parallel each other. 
Note $\partial N(Q_i)\cong S^1\x S^2$. 
We can suppose the following. 

\noindent 
(1) 
The intersection (each $E_{Yi}$) $\cap\partial N(Q_i)$ is 
a 2-spehre. 

\noindent 
(2) 
This 2-sphere is 
 $p\x S^2\subset S^1\x S^2\cong\partial N(Q_i)$, 
 where 
 $p$ is a point $\in$ (this $S^1$). 
 
\noindent 
(3)
Each 2-sphere is parallel to other 2-spheres in 
$S^1\x S^2\cong\partial N(Q_i)$. 

\vskip3mm

We give an orientation to each $E_{Yi}$ so that 
the orientation of 
$P_i\cap E_{Yi}$
is the opposite one of 
that of $P_i\cap M_{*i}$. 
Note that we do not suppose the orientation of 
$E_{Yi}$ coincides with that of 
$Q_i$.

Carry out surgeries on \newline
`$M_{*i} \amalg (\nu$ copies of $E_{Yi})$' 
by using 4-dimensional $\nu$ 1-handles with the following properties. 
(Note that $M_{*i}\cap E_{Yi}=\phi$ by the above Claim.)

\noindent(1) 
The handles are embedded in 
$Y_i$.

\noindent(2) 
The handles are attached along two points 
such that one point is in 
$P_i\cap ({\text{$\nu$ copies of }} E_{Yi})$ 
and that the other point is in 
$P_i\cap M_{*i}$. 
 (Note that the orientation of $M_{*i}$ is compatible  with that of $E_{Yi}$.)

\noindent(3) 
The core of each 1-handle is 
a 1-dimensional 
`submanifold with boundary' of $P_i$. 
(Note $P_i\cong S^1$.)



\vskip3mm

Thus we made a new submanifold 
$M^\#_{*i}$ from 

\noindent 
$M_{*i}$, `$\nu$ copies of $E_{Yi}$', and `the above $\nu$ 1-handles'
 
\noindent by these surgeries. 
Note that these surgeries can avoid making any self-intersection 
of $M^\#_{*i}$ for each $i$ by using the above 1-hndles appropriately. 
Note that 
$\partial M^\#_{*i}-\partial Y_i=
 \partial E_{Yi}-\partial Y_i$    
and that 
$\partial M^\#_{*i}-\partial Y_i$ 
is a set of the 2-spheres $\subset \partial N(Q_i)$.

\vskip3mm
\noindent 
Claim.{\it These surgeries do not change 
$I\subset M_{\alpha i}\cap M_{\beta i}\cap M_{\gamma i}$ 
for each $\{\alpha, \beta, \gamma\}$.}
\vskip3mm

\noindent 
Proof.
If necessary, we carry out these surgeries on all $M_{*i}(i=1,...,m)$.  
These surgeries are done 
in the interior of the compact set $N(P_i)$, 
where $N(P_i)$ is the tubular neighborhood of $P_i$ in $Y_i$. 
Suppose that  the triple point set 
$\widetilde{M_\alpha}\cap \widetilde{M_\beta}\cap \widetilde{M_\gamma}$ 
intersects $N(P_i)$ before these surgeries. 
By using isotopy of $\widetilde{M_*}$, 
we  can suppose that 
$\widetilde{M_\alpha}\cap \widetilde{M_\beta}\cap \widetilde{M_\gamma}$
does not intersect $N(P_i)$. 
Suppose that 
the triple point set 
$\widetilde{M_\alpha}\cap \widetilde{M_\beta}\cap \widetilde{M_\gamma}$ 
is obtained after these surgeries. 
Then 
the new triple point set 
$\widetilde{M_\alpha}\cap \widetilde{M_\beta}\cap \widetilde{M_\gamma}$ 
is a disjoint union of circles. 
Because we have the following. 

\noindent 
(1) 
These surgeries are done in the interior of the compact set 
$N(P_i)$.

\noindent 
(2)
The triple point set does not exist in $N(P_i)$ before these surgeries. 
(Because of the above procedure.)  

\vskip3mm
This completes the proof. 
\vskip3mm

Carry out the ribbon-move surgery along $P_i$ and $Q_i$.
Then, by Claim \ref{coveringsurgery}, we have the following. 

\noindent 
(1) 
$\widetilde{X^\infty_{K}}$ 
is changed into 
$\widetilde{X^\infty_{K'}}$.

\noindent 
(2) 
Let $Y_{Ki}$ be changed into $Y'_{Ki}$. 
Then 
$Y'_{Ki}$ is $Y_{K'i}$. 
Recall that $Y_{Ki}$ and $Y_i$ are same 
(it is written in \S\ref{ccc} and before Definition \ref{subsur}.). 

\vskip3mm

After this ribbon-move surgery,
 we can carry out the following surgeries.

 
Take 3-dimensional 3-handles in new $Y'_{Ki}=Y_{K'i}$. 
Strictly to say, the 3-dimensional 3-handles are in 
$Y'_{Ki}-(Y'_{Ki}\cap Y_{Ki})$. 
Attach these 3-dimensional 3-handles to $M^\#_{*i}$ 
along all of the 2-spheres 
$\partial M^\#_{*i}-\partial Y_{Ki}=
 \partial E_{Yi}-\partial Y_{Ki}$. 
Thus we obtain a new 3-manifold $M'_{*i}$ with boundary.

When we change $M_{*i}$, we do not change \newline
$I (\subset \widetilde{M_\alpha}\cap 
\widetilde{M_\beta}\cap 
\widetilde{M_\gamma}$ ),  
$\partial M_{*i} (\subset \partial Y_{Ki})$, 
or
$(\partial E_{Y_i}\cap \partial Y_{Ki})$.  

Note: In the above procedure, we may move $I$ by isotopy. 
However, we do not change the diffeomorphism type of $I$. 
Furthermore, we do not move $\partial I$.

Hence we have the following. 

\noindent 
(1) 
$\partial M'_{*i}=\partial M_{*i}\amalg(\partial E_{Y_i}\cap \partial Y_{Ki})$.  

\noindent 
(2) 
$\widetilde{M_*}'
   =(\amalg_{i=-\infty}^{\infty} M'_{*i})
\cup(\amalg_{i=-\infty}^{\infty} 
\{\nu\mathrm{\hskip1mm copies\hskip1mm of\hskip1mm} E_{NVi}\})
\cup(\amalg_{i=-\infty}^{\infty} M_{*NVi})
$ 

\noindent
is an open 3-manifold without boundary. 
(Note that $E_{NVi}$ is defined 
after a few lines from the proof of Claim \ref{E}.
Note that $M_{*NVi}$ is defined 
before a few lines from Claim \ref{h1saketa}. )

$\widetilde{M_*}' $ is a submanifold of $\widetilde{X^\infty_{K'}}$. 
$\widetilde{M_*}'$ 
represents an element 
$\in H^\infty_3(\widetilde{X^\infty_{K'}};\R)$. 


\noindent 
(3) 
$\{\widetilde{M_1}',...,\widetilde{M_m}'\} (m\in\N)$ 
is an $o$-set for $K'$.

\noindent 
(4)
For each set $\{\alpha, \beta, \gamma\}\subset \{1,...,m\}$, 
we have 

\noindent 
$o(\widetilde{M_\alpha}', \widetilde{M_\beta}', \widetilde{M_\gamma}')$= 
$o(\widetilde{M_\alpha}, \widetilde{M_\beta}, \widetilde{M_\gamma})$.














\section{Problems}\label{problems}

\noindent
Here, we submit Problem  \ref{problem} again.

\vskip3mm
\begin{prob}\label{problem1} (essentially same as Problem \ref{problem}.)   
Classify 2-knots by the ribbon-move equivalence.
\end{prob}
\vskip3mm

In particular, the following problems interest us. 

\vskip3mm
\begin{prob}\label{problem2}
Is there a nontrivial 2-knot $K$ with the following properties? 

\noindent(1) 
A Seifert hypersurface of $K$ is a punctured integral homology sphere. 

\noindent(2) 
$K$ is ribbon-move equivalent to the trivial knot. 
\end{prob}
\vskip3mm

\begin{prob}\label{problem3}
Is there a nontrivial 2-knot $K$ with the following properties?

\noindent(1) 
A Seifert hypersurface of $K$ is a punctured integral homology sphere. 

\noindent(2) 
$\mu(K)$=0. 

\noindent(3) 
$K$ is not ribbon-move equivalent to the trivial knot. 

\end{prob}

\vskip3mm
\cite{Zeeman} proved that 
a Seifert hypersurface of the five twist spun knot of the trefoil knot 
is the punctured Poincar\'e homology sphere. 
(Note \ref{Z} quotes this result.) 
This 2-knot is called $Z$.

\vskip3mm
\begin{prob}\label{problem4}
Is  
$
\left\{
\begin{array}{l}
 Z\# Z \\
 Z\#(-Z)
 \end{array}
\right.
$ 
ribbon-move equivalent to the trivial knot?
(Note that $-Z$ is the 2-knot which has the opposite orientation of $Z$.) 
\end{prob}
\vskip3mm



We introduce another local move.

\vskip3mm
\begin{defn}\label{XO} 
Let $K$ be a 2-knot $\subset S^4$. 
Embed $S^1\x D^3$ trivially in $S^4$, 
where 
$S^1$ is a circle and where $D^3$ is a close 3-disc. 
Suppose that the following hold. 

\noindent (1) 
$K\cap (S^1\x D^3)$ is 
$(S^1\x I)\amalg(S^1\x I)$, where $I$ is the interval. 

\noindent (2) 
$K\cap (S^1\x D^3)$ is

\vskip3mm
\unitlength 0.1in
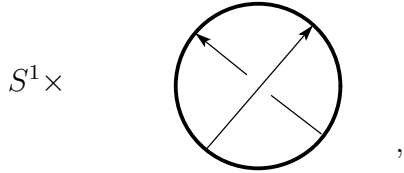
\begin{picture}(20.30,8.62)(13.70,-15.61)
\put(13.7000,-11.7000){\makebox(0,0)[lb]{$S^1\x$}}%
%
\special{pn 20}%
\special{ar 2680 1130 431 431  0.8264709 6.2831853}%
\special{ar 2680 1130 431 431  0.0000000 0.8017901}%
%
\special{pn 8}%
\special{pa 2750 1180}%
\special{pa 3010 1380}%
\special{fp}%
%
\special{pn 8}%
\special{pa 2620 1070}%
\special{pa 2360 860}%
\special{fp}%
\special{sh 1}%
\special{pa 2360 860}%
\special{pa 2399 917}%
\special{pa 2401 894}%
\special{pa 2424 886}%
\special{pa 2360 860}%
\special{fp}%
%
\special{pn 8}%
\special{pa 2410 1460}%
\special{pa 2960 820}%
\special{fp}%
\special{sh 1}%
\special{pa 2960 820}%
\special{pa 2901 858}%
\special{pa 2925 860}%
\special{pa 2932 884}%
\special{pa 2960 820}%
\special{fp}%
\put(34.0000,-15.1000){\makebox(0,0)[lb]{,}}%
\end{picture}%
\vskip3mm

\noindent 
where we have the following.

\noindent (i) 
The bold line and its interior in the above figure 
represent the 3-disc $D^3$ 

\noindent (ii) 
The arrows of finite lines represent 
a submanifold of $K$. 

\noindent (iii) 
$S^1\x$ (each of the two arrows) means  
each of the above $(S^1\x I)\amalg(S^1\x I)$. 

\vskip3mm
\noindent
Fix this chart of $S^1\x D^3\subset \R^4$. 

Let  $K'$ be a 2-knot with the following properties. 

\noindent(1)
$K\cap (S^1\x D^3)$ is 

\vskip3mm
\unitlength 0.1in
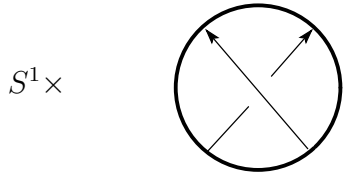
\begin{picture}(17.41,8.62)(13.70,-15.61)
\put(13.7000,-11.7000){\makebox(0,0)[lb]{$S^1\x$}}%
%
\special{pn 20}%
\special{ar 2680 1130 431 431  0.8264709 6.2831853}%
\special{ar 2680 1130 431 431  0.0000000 0.8017901}%
%
\special{pn 8}%
\special{pa 2940 1450}%
\special{pa 2410 830}%
\special{fp}%
\special{sh 1}%
\special{pa 2410 830}%
\special{pa 2438 894}%
\special{pa 2445 871}%
\special{pa 2469 868}%
\special{pa 2410 830}%
\special{fp}%
%
\special{pn 8}%
\special{pa 2410 1470}%
\special{pa 2630 1230}%
\special{fp}%
%
\special{pn 8}%
\special{pa 2750 1070}%
\special{pa 2960 830}%
\special{fp}%
\special{sh 1}%
\special{pa 2960 830}%
\special{pa 2901 867}%
\special{pa 2925 870}%
\special{pa 2931 893}%
\special{pa 2960 830}%
\special{fp}%
\end{picture}%
\vskip3mm

\noindent
in the above chart. 

\noindent(2) 
$K'\cap\overline{S^4-(S^1\x D^3)}$
=
$K\cap\overline{S^4-(S^1\x D^3)}$. 

Then we say that 
$K'$ is obtained from $K$ by one {\it XO-move}.

If $K''$ is obtained from $K$ by a sequence of XO-moves 
then we say that 
$K''$ is {\it XO-move equivalent} to $K$. 
\end{defn}

\vskip3mm 

All $n$-twist spun knots could be XO-move equivalent to the trivial knot. 
The obstructions for ribbon-moves in \cite{O1,O2} could not be obstructions 
for XO-moves.

\vskip3mm
\begin{prob}\label{XOQ}
\noindent(1) 
Are all 2-knots XO-move equivalent to the trivial knot?

\noindent(2)
Is the $o$-invariant an obstruction for XO-moves?
\end{prob}
\vskip3mm



\end{document}